# A TEST FOR MODEL SPECIFICATION OF DIFFUSION PROCESSES

By Song Xi Chen,[1,2] Jiti Gao[1,3] and Cheng Yong Tang[1,2]

*Iowa State University, University of Western Australia
and Iowa State University*

We propose a test for model specification of a parametric diffusion process based on a kernel estimation of the transitional density of the process. The empirical likelihood is used to formulate a statistic, for each kernel smoothing bandwidth, which is effectively a Studentized $L_2$-distance between the kernel transitional density estimator and the parametric transitional density implied by the parametric process. To reduce the sensitivity of the test on smoothing bandwidth choice, the final test statistic is constructed by combining the empirical likelihood statistics over a set of smoothing bandwidths. To better capture the finite sample distribution of the test statistic and data dependence, the critical value of the test is obtained by a parametric bootstrap procedure. Properties of the test are evaluated asymptotically and numerically by simulation and by a real data example.

**1. Introduction.** Let $X_1, \ldots, X_{n+1}$ be $n+1$ equally spaced (with spacing $\Delta$ in time) observations of a diffusion process

$$(1.1) \qquad dX_t = \mu(X_t)\,dt + \sigma(X_t)\,dB_t,$$

where $\mu(\cdot)$ and $\sigma^2(\cdot) > 0$ are, respectively, the drift and diffusion functions, and $B_t$ is the standard Brownian motion. Suppose a parametric specification of model (1.1) is

$$(1.2) \qquad dX_t = \mu(X_t; \theta)\,dt + \sigma(X_t; \theta)\,dB_t,$$

Received November 2005; revised February 2007.

[1]Supported by the National University of Singapore Academic research grant.

[2]Supported by the NSF Grants DMS-06-04563 and SES-05-18904.

[3]Supported by the Australian Research Council Discovery Grant.

where $\theta$ is a parameter within a parameter space $\Theta \subset R^d$ for a positive integer $d$. The focus of this paper is on testing the validity of the parametric specification (1.2) based on a set of discretely observed data $\{X_{t\Delta}\}_{t=1}^{n+1}$.

In a pioneer work that represents a break-through in financial econometrics, Aït-Sahalia [1] considered two approaches for testing the parametric specification (1.2). The first one was based on a $L_2$-distance between a kernel stationary density estimator and the parametric stationary density implied by model (1.2) with the critical value of the test obtained from the asymptotic normal distribution of the test statistic. The advantage of the test is that the parametric stationary density is easily derivable for almost all processes and performing the test is straightforward. There are several limitations with the test. One is, as pointed out by Aït-Sahalia [1], that a test that targets on the stationary distribution is not conclusive, as different processes may share a common stationary distribution. Another is that it can take a long time for a process is produce a sample path that contains enough information for accurate estimation of the stationary distribution. These were confirmed by Pritsker [42] who reported noticeable discrepancy between the simulated and nominal sizes of the test under a set of Vasicek [47] diffusion processes. In the same paper Aït-Sahalia considered another approach based on certain discrepancy measure regarding the transitional distribution of the process derived from the Kolmogorov-backward and backward equations. The key advantage of a test that targets on the transitional density is that it is conclusive as transitional density fully specifies the dynamics of a diffusion process due to its Markovian property.

In this paper we propose a test that is focused on the specification of the transitional density of a process. The basic building blocks used in constructing the test statistic are the kernel estimator of the transitional density function and the empirical likelihood (Owen [41]). We first formulate an integrated empirical likelihood ratio statistic for each smoothing bandwidth used in the kernel estimator, which is effectively a $L_2$-distance between the kernel transitional density estimator and the parametric transitional density implied by the process. The use of the empirical likelihood allows the $L_2$-distance being standardized by the variation. We also implement a series of measures to make the test work more efficiently. This includes properly smoothing the parametric transitional density so as to cancel out the bias induced by the kernel estimation, which avoids undersmoothing and simplifies theoretical analysis. To make the test robust against the choice of smoothing bandwidth, the test statistic is formulated based on a set of bandwidths. Finally, a parametric bootstrap procedure is employed to obtain the critical value of the test, and to better capture the finite sample distribution of the test statistic and data dependence induced by the stochastic process.

A continuous-time diffusion process and discrete-time time series share some important features. They can be both Markovian and weakly dependent satisfying certain mixing condition. The test proposed in this paper draws experiences from research works on kernel based testing and estimation of discrete time models established in the last decade or so. Kernel-based tests have been shown to be effective in testing discrete time series models as demonstrated in Robinson [44], Fan and Li [24], Hjellvik, Yao and Tjøstheim [30], Li [38], Aït-Sahalia, Bickel and Stoker [4], Gozalo and Linton [27] and Chen, Härdle and Li [12]; see Hart [29] and Fan and Yao [18] for extended reviews and lists of references. For kernel estimation of diffusion processes, in addition to Aït-Sahalia [1], Jiang and Knight [35] proposed a semiparametric kernel estimator; Fan and Zhang [21] examined the effects of high order stochastic expansions and proposed separate generalized likelihood ratio tests for the drift and diffusion functions; Bandi and Phillips [7] considered a two-stage kernel estimation without the strictly stationary assumption. See Cai and Hong [9] and Fan [15] for comprehensive reviews.

In an important development after Aït-Sahalia [1], Hong and Li [32] developed a test for diffusion processes via a conditional probability integral transformation. The test statistic is based on a $L_2$-distance between the kernel density estimator of the transformed data and the uniform density implied under the hypothesized model. Although the kernel estimator is employed, the transformation leads to asymptotically independent uniform random variables under the hypothesized model. Hence, the issue of modeling data dependence induced by diffusion processes is avoided. In a recent important development, Aït-Sahalia, Fan and Peng [5] proposed a test for the transitional densities of diffusion and jump diffusion processes based on a generalized likelihood ratio statistic, which, like our current proposal, is able to fully test diffusion process specification and has attractive power properties.

The paper is structured as follows. Section 2 outlines the hypotheses and the kernel smoothing of transitional densities. The proposed EL test is given in Section 3. Section 4 reports the main results of the test. Section 5 considers computational issues. Results from simulation studies are reported in Section 6. A Federal fund rate data set is analyzed in Section 7. All technical details are given in the Appendix.

## 2. The hypotheses and kernel estimators.

Let $\pi(x)$ be the stationary density and $p(y|x; \Delta)$ be the transitional density of $X_{(t+1)\Delta} = y$ given $X_{t\Delta} = x$ under model (1.1), respectively; and $\pi_\theta(x)$ and $p_\theta(y|x, \Delta)$ be their parametric counterparts under model (1.2). To simplify notation, we suppress $\Delta$ in the notation of transitional densities and write $\{X_{t\Delta}\}$ as $\{X_t\}$ for the observed data. Let $\mathcal{X}$ be the state space of the process.

Although $\pi_\theta(x)$ has a close form expression via Kolmogorov forward equation

$$\pi_\theta(x) = \frac{\xi(\theta)}{\sigma^2(x, \theta)} \exp\left\{\int_{x_0}^x \frac{2\mu(t, \theta)}{\sigma^2(t, \theta)} \, dt\right\},$$

where $\xi(\theta)$ is a normalizing constant, $p_\theta(y|x)$ as defined by the Kolmogorov-backward equation may not admit a close form expression. However, this problem is overcome by Edgeworth type approximations developed by Aït-Sahalia [2, 3]. As the transitional density fully describes the dynamics of a diffusion process, the hypotheses we would like to test are

$$H_0 : p(y|x) = p_{\theta_0}(y|x) \text{ for some } \theta_0 \in \Theta \text{ and all } (x,y) \in S \subset \mathcal{X}^2 \quad \text{versus}$$

$$H_1 : p(y|x) \neq p_\theta(y|x) \text{ for all } \theta \in \Theta \text{ and some } (x,y) \in S \subset \mathcal{X}^2,$$

where $S$ is a compact set within $\mathcal{X}^2$ and can be chosen based on the kernel transitional density estimator given in (2.1) below; see also demonstrations in simulation and case studies in Sections 6 and 7. As we are to properly smooth the parametric density $p_\theta(y|x)$, the boundary bias associated with the kernel estimators ([16] and [40]) is avoided.

Let $K(\cdot)$ be a kernel function which is a symmetric probability density function, $h$ be a smoothing bandwidth such that $h \to 0$ and $nh^2 \to \infty$ as $n \to \infty$, and $K_h(\cdot) = h^{-1}K(\cdot/h)$. The kernel estimator of $p(y|x)$ is

$$(2.1) \qquad \hat{p}(y|x) = n^{-1} \sum_{t=1}^{n} K_h(x - X_t) K_h(y - X_{t+1}) / \hat{\pi}(x),$$

where $\hat{\pi}(x) = (n+1)^{-1} \sum_{t=1}^{n+1} K_h(x - X_t)$ is the kernel estimator of the stationary density used in Aït-Sahalia [1]. The local polynomial estimator introduced by Fan, Yao and Tong [19] can also be employed without altering the main results of this paper. It is known (Hydman and Yao [34]) that

$$E\{\hat{p}(y|x) - p(y|x)\} = \frac{1}{2}\sigma_k^2 h^2 \left( \frac{\partial^2 p(y|x)}{\partial x^2} + \frac{\partial^2 p(y|x)}{\partial y^2} + 2\frac{\pi'(x)}{\pi(x)} \frac{\partial p(y|x)}{\partial x} \right)$$
$$+ o(h^2),$$

$$\mathrm{Var}\{\hat{p}(y|x)\} = \frac{R^2(K)p(y|x)}{nh^2\pi(x)}(1 + o(1)),$$

where $\sigma_K^2 = \int u^2 K(u)\,du$ and $R(K) = \int K^2(u)\,du$. Here we use a single bandwidth $h$ to smooth the bivariate data $(X_t, X_{t+1})$. This is based on a consideration that both $X_t$ and $X_{t+1}$ are identically distributed and hence have the same scale which allows one smoothing bandwidth to smooth for both components. Nevertheless, the results in this paper can be generalized to the situation where two different bandwidths are employed.

Let $\tilde{\theta}$ be a consistent estimator of $\theta$ under model (1.2), for instance, the maximum likelihood estimator under $H_0$, and

$$(2.2) \qquad w_t(x) = K_h(x - X_t) \frac{s_{2h}(x) - s_{1h}(x)(x - X_t)}{s_{2h}(x)s_{0h}(x) - s_{1h}^2(x)}$$

be the local linear weight with $s_{rh}(x) = \sum_{s=1}^{n} K_h(x - X_s)(x - X_s)^r$ for $r = 0, 1$ and 2. In order to cancel out the bias in $\hat{p}(y|x)$, we smooth $p_{\hat{\theta}}(y|x)$ as

$$(2.3) \qquad \tilde{p}_{\hat{\theta}}(y|x) = \frac{\sum_{t=1}^{n+1} K_h(x - X_t) \sum_{s=1}^{n+1} w_s(y) p_{\hat{\theta}}(X_s|X_t)}{\sum_{t=1}^{n+1} K_h(x - X_t)}.$$

Here we apply the kernel smoothing twice: first for each $X_t$ using the local linear weight to smooth $p_{\hat{\theta}}(X_s|X_t)$ and then employing the standard kernel to smooth with respect to $X_t$. This is motivated by Härdle and Mammen [28]. It can be shown from the standard derivations in Fan and Gijbels [16] that, under $H_0$,

$$(2.4) \qquad E\{\hat{p}(y|x) - \tilde{p}_{\hat{\theta}}(y|x)\} = o(h^2)$$

and

$$(2.5) \qquad \mathrm{Var}\{\hat{p}(y|x) - \tilde{p}_{\hat{\theta}}(y|x)\} = \mathrm{Var}\{\hat{p}(y|x)\}\{1 + o(1)\}.$$

Hence, the biases of $\hat{p}(y|x)$ and $\tilde{p}_{\hat{\theta}}(y|x)$ are canceled out in the leading order, while smoothing the parametric density does not affect the asymptotic variance.

## 3. Formulation of test statistic.

The test statistic is formulated by the empirical likelihood (EL) (Owen [41]). Despite its being intrinsically non-parametric, EL possesses two key properties of a parametric likelihood: the Wilks' theorem and the Bartlett correction. Qin and Lawless [43] established EL for parameters defined by generalized estimating equations which is the broadest framework for EL formulation so far, which was extended by Kitamura [36] to dependent observations. Chen and Cui [10] showed that the EL admits Bartlett correction under this general framework. See also Hjort, McKeague and Van Keilegom [31] for extensions. The EL has been used for goodness-of-fit tests of various model structures. Fan and Zhang [23] proposed a sieve EL test for a varying-coefficient regression model that extends the test of Fan, Zhang and Zhang [22]; Tripathi and Kitamura [46] studied a test for conditional moment restrictions; Chen, Härdle and Li [12] proposed an EL test for time series regression models. See also [37] for survival data.

We now formulate the EL for the transitional density at a fixed $(x, y)$. For $t = 1, \ldots, n$, let $q_t(x, y)$ be nonnegative weights allocated to $(X_t, X_{t+1})$. The EL evaluated at $\tilde{p}_{\hat{\theta}}(y|x)$ is

$$(3.1) \qquad L\{\hat{p}_{\hat{\theta}}(y|x)\} = \max \prod_{t=1}^{n} q_t(x, y)$$

subject to $\sum_{t=1}^{n} q_t(x, y) = 1$ and

$$(3.2) \qquad \sum_{t=1}^{n} q_t(x, y) K_h(x - X_t) K_h(y - X_{t+1}) = \tilde{p}_{\hat{\theta}}(y|x)\hat{\pi}(x).$$

By introducing a Lagrange multiplier $\lambda(x, y)$, the optimal weights as solutions to (3.1) and (3.2) are

$$(3.3) \qquad q_t(x, y) = n^{-1}\{1 + \lambda(x, y)T_t(x, y)\}^{-1},$$

where $T_t(x, y) = K_h(x - X_t)K_h(y - X_{t+1}) - \tilde{p}_{\hat{\theta}}(x, y)$ and $\lambda(x, y)$ is the root of

$$(3.4) \qquad \sum_{t=1}^{n} \frac{T_t(x, y)}{1 + \lambda(x, y)T_t(x, y)} = 0.$$

The overall maximum EL is achieved at $q_t(x, y) = n^{-1}$ which maximizes (3.1) without constraint (3.2). Hence, the log-EL ratio is

$$(3.5) \qquad \begin{aligned} \ell\{\tilde{p}_{\hat{\theta}}(y|x)\} &= -2\log([L\{\tilde{p}_{\hat{\theta}}(y|x)\}n^n]) \\ &= 2\sum \log\{1 + \lambda(x, y)T_t(x, y)\}. \end{aligned}$$

It may be shown by similar derivations to those in Chen, Härdle and Li [12] that

$$(3.6) \qquad \sup_{(x, y) \in S} |\lambda(x, y)| = o_p\{(nh^2)^{-1/2}\log(n)\}.$$

Let $\bar{U}_1(x, y) = (nh^2)^{-1}\sum T_t(x, y)$ and $\bar{U}_2(x, y) = (nh^2)^{-1}\sum T_\ell^2(x, y)$. From (3.4) and (3.6), $\lambda(x, y) = \bar{U}_1(x, y)\bar{U}_2^{-1}(x, y) + O_p\{(nh^2)^{-1}\log^2(n)\}$ uniformly with respect to $(x, y) \in S$. This leads to

$$(3.7) \qquad \begin{aligned} \ell\{\tilde{p}_{\hat{\theta}}(y|x)\} &= nh^2\bar{U}_1^2(x, y)\bar{U}_2^{-1}(x, y) + O_p\{n^{-1/2}h^{-1/2}\log^3(n)\} \\ &= nh^2\frac{\{\hat{p}(y|x) - \tilde{p}_{\hat{\theta}}(y|x)\}^2}{V(y|x)} + O_p\{h^2 + n^{-1/2}h^{-1/2}\log^3(n)\} \end{aligned}$$

uniformly for $(x, y) \in S$, where $V(y|x) = R^2(K)p(y|x)\pi^{-1}(x)$. Hence, the EL ratio is a Studentized local goodness-of-fit measure between $\hat{p}(y|x)$ and $\tilde{p}_{\hat{\theta}}(y|x)$ as $\text{Var}\{\hat{p}(y|x)\} = (nh^2)^{-1}V(y|x)$.

Integrating the EL ratio against a weight function $\omega(\cdot, \cdot)$ supported on $S$, the global goodness-of-fit measure based on a single bandwidth is

$$(3.8) \qquad N(h) = \iint \ell\{\tilde{p}_{\hat{\theta}}(y|x)\}\omega(x, y)\,dx\,dy.$$

To make the test less dependent on a single bandwidth $h$, we compute $N(h)$ over a bandwidth set $\mathcal{H} = \{h_k\}_{k=1}^{J}$, where $h_k/h_{k+1} = a$ for some $a \in (0, 1)$. The choice of $\mathcal{H}$ can be guided by the cross-validation method of Fan and Yim [20] or other bandwidth selection methods; see Section 5 for more discussions and demonstration. This formulation is motivated by Fan [14] and Horowitz and Spokoiny [33], both considered achieving the optimal convergence rate for the distance between a null hypothesis and a series of local alternative hypotheses in testing regression models. To our best knowledge, Fan [14] was the first to propose the adaptive test and showed its oracle property. The adaptive result is more explicitly given in Fan and Huang [17]. Fan, Zhang and Zhang [22] also explicitly adapted the multi-frequency test of Fan [14] into the multi-scale test and obtained the adaptive minimax result. As we are concerned with testing against a fixed alternative only, it is adequate to have a finite number of bandwidths in $\mathcal{H}$ in our context.

The final test statistic based on the bandwidth set $\mathcal{H}$ is

$$(3.9) \qquad L_n = \max_{1 \le k \le J} \frac{N(h_k) - 1}{\sqrt{2} h_k},$$

where the standardization reflects that $\mathrm{Var}\{N(h)\} = O(2h^2)$ as shown in the [Appendix](#).

## 4. Main results.
Our theoretical results are based on the following assumptions.

ASSUMPTION 1. (i) The process $\{X_t\}$ is strictly stationary and $\alpha$-mixing with mixing coefficient $\alpha(t) \le C_\alpha \alpha^t$, where

$$\alpha(t) = \sup\{|P(A \cap B) - P(A)P(B)| : A \in \Omega_1^s, B \in \Omega_{s+t}^\infty\}$$

for all $s, t \ge 1$, where $C_\alpha$ is a finite positive constant, $\Omega_i^j$ denotes the $\sigma$-field generated by $\{X_t : i \le t \le j\}$, and $\alpha$ is a constant in $(0, 1)$.

(ii) $K(\cdot)$ is a bounded symmetric probability density supported on $[-1, 1]$ and has bounded second derivative; and $\omega(x, y)$ is a bounded probability density supported on $S$.

(iii) For the bandwidth set $\mathcal{H}$, $h_1 = c_1 n^{-\gamma_1}$ and $h_J = c_J n^{-\gamma_2}$, in which $\frac{1}{7} < \gamma_2 \le \gamma_1 < \frac{1}{4}$, $c_1$ and $c_J$ are constants satisfying $0 < c_1, c_J < \infty$, and $J$ is a positive integer not depending on $n$.

ASSUMPTION 2. (i) Each of the diffusion processes given in [(1.1)](#) and [(1.2)](#) admits a unique weak solution and possesses a transitional density with $p(y|x) = p(y|x, \Delta)$ for [(1.1)](#) and $p_\theta(y|x) = p_\theta(y|x, \Delta)$ for [(1.2)](#).

(ii) Let $p_{s_1, s_2, \ldots, s_l}(\cdot)$ be the joint probability density of $(X_{1+s_1}, \ldots, X_{1+s_l})$. Assume that each $p_{s_1, \tau_2, \ldots, s_l}(x)$ is three times differentiable in $x \in \mathcal{X}^l$ for $1 \le l \le 6$.

(iii) The parameter space $\Theta$ is an open subset of $R^d$ and $p_\theta(y|x)$ is three times differentiable in $\theta \in \Theta$. For every $\theta \in \Theta$, $\mu(x; \theta)$ and $\sigma^2(x; \theta)$, and $\mu(x)$ and $\sigma^2(x)$ are all three times continuously differentiable in $x \in \mathcal{X}$, and both $\sigma(x)$ and $\sigma(x; \theta)$ are positive for $x \in S$ and $\theta \in \Theta$.

ASSUMPTION 3. (i) $E[(\frac{\partial p_\theta(X_{t+1}|X_t)}{\partial \theta})(\frac{\partial p_\theta(X_{t+1}|X_t)}{\partial \theta})^\tau]$ is of full rank. Let $G(x, y)$ be a positive and integrable function with $E[\max_{1 \le t \le n} G(X_t, X_{t+1})] < \infty$ uniformly in $n \ge 1$ such that $\sup_{\theta \in \Theta} |p_\theta(y|x)|^2 \le G(x, y)$ and $\sup_{\theta \in \Theta} \| \triangledown_\theta^j p_\theta(y|x) \|^2 \le G(x, y)$ for all $(x, y) \in S$ and $j = 1, 2, 3$, where $\triangledown_\theta p_\theta(\cdot|\cdot) = \frac{\partial p_\theta(\cdot|\cdot)}{\partial \theta}$, $\triangledown_\theta^2 p_\theta(\cdot|\cdot) = \frac{\partial^2 p_\theta(\cdot|\cdot)}{(\partial \theta)^2}$ and $\triangledown_\theta^3 p_\theta(\cdot|\cdot) = \frac{\partial^3 p_\theta(\cdot|\cdot)}{(\partial \theta)^3}$.

(ii) $p(y|x) > c_1 > 0$ for all $(x, y) \in S$ and the stationary density $\pi(x) > c_2 > 0$ for all $x \in S_x$ which is the projection of $S$ on $\mathcal{X}$.

ASSUMPTION 4. Under either $H_0$ or $H_1$, there is a $\theta^* \in \Theta$ and a sequence of positive constants $\{a_n\}$ that diverges to infinity such that, for any $\varepsilon > 0$ and some $C > 0$, $\lim_{n \to \infty} P(a_n \|\tilde{\theta} - \theta^*\| > C) < \varepsilon$ and $\sqrt{n} h a_n^{-1} = o(\sqrt{h})$ for any $h \in \mathcal{H}$.

Assumption 1(i) imposes the strict stationarity and $\alpha$-mixing condition on $\{X_t\}$. Under certain conditions, such as Assumption A2 of Aït-Sahalia [1] and Conditions (A4) and (A5) of Genon-Catalot, Jeantheau and Larédo [26], Assumption 1(i) holds. Assumption 1(ii) and (iii) are quite standard conditions imposed on the kernel and the bandwidth in kernel estimation. Assumption 2 is needed to ensure the existence and uniqueness of a solution and the transitional density function of the diffusion process. Such an assumption may be implied under Assumptions 1–3 of Aït-Sahalia [3], which also cover nonstationary cases. For the stationary case, Assumptions A0 and A1 of Aït-Sahalia [1] ensure the existence and uniqueness of a stationary solution of the diffusion process. Assumption 3 imposes additional conditions to ensure the smoothness of the transitional density and the identifiability of the parametric transitional density. The $\theta^*$ in Assumption 4 is the true parameter $\theta_0$ under $H_0$. When $H_1$ is true, $\theta^*$ can be regarded as a projection of the parameter estimator $\tilde{\theta}$ onto the null parameter space. Assumption 4 also requires that $a_n$, the rate convergence of $\tilde{\theta}$ to $\theta^*$, is faster than $\sqrt{n} h$, the convergence rate for the kernel transitional density estimation. This is certainly satisfied when $\tilde{\theta}$ converges at the rate of $\sqrt{n}$ as attained by the maximum likelihood estimation. Our use of the general convergence rate $a_n$ for the parameter estimation is to cover situations where the parameter estimator has a slower rate than $\sqrt{n}$, for instance, when estimation is based on certain forms of discretization which requires $\Delta \to 0$ in order to be consistent (Lo [39]).

Let $K^{(2)}(z, c) = \int K(u) K(z + cu) \, du$, a generalization to the convolution of $K$, $\nu(t) = \int \{K^{(2)}(tu, t)\}^2 \, du \int \{K^{(2)}(v, t)\}^2 \, dv$ and

$$\Sigma_J = \frac{2}{R^4(K)} \int\int \omega^2(x, y) \, dx \, dy \, (\nu(a^{i-j}))_{J \times J}$$

be a $J \times J$ matrix, where $a$ is the fixed factor used in the construction of $\mathcal{H}$. Furthermore, let $\mathbf{1}_J$ be a $J$-dimensional vector of ones and $\beta = \frac{1}{R(K)} \iint \frac{p(x,y)}{\pi(y)} \omega(x, y) \, dx \, dy$.

THEOREM 1. *Under Assumptions 1–4 and $H_0$, $L_n \xrightarrow{d} \max_{1 \leq k \leq J} Z_k$ as $n \to \infty$ where $Z = (Z_1, \ldots, Z_J)^T \sim N(\beta \mathbf{1}_J, \Sigma_J)$.*

Theorem 1 brings a little surprise in that the mean of $Z$ is nonzero. This is because, although the variance of $\tilde{p}_\theta(x, y)$ is at a smaller order than that of $\hat{p}(x, y)$, it contributes to the second-order mean of $N(h)$ which emerges after dividing $\sqrt{2}h$ in (3.9). However, this does not affect $L_n$ being a test statistic.

We are reluctant to formulate a test based on Theorem 1 as the convergence would be slow. Instead, we propose the following parametric bootstrap procedure to approximate $l_\alpha$, the $1 - \alpha$ quantile of $L_n$ for a nominal significance level $\alpha \in (0, 1)$:

STEP 1. Generate an initial value $X_0^*$ from the estimated stationary density $\pi_{\tilde{\theta}}(\cdot)$. Then simulate a sample path $\{X_t^*\}_{t=1}^{n+1}$ at the same sampling interval $\Delta$ according to $dX_t = \mu(X_t; \tilde{\theta}) \, dt + \sigma(X_t; \tilde{\theta}) \, dB_t$.

STEP 2. Let $\tilde{\theta}^*$ be the estimate of $\theta$ based on $\{X_t^*\}_{t=1}^{n+1}$. Compute the test statistic $L_n$ based on the resampled path and denote it by $L_n^*$.

STEP 3. For a large positive integer B, repeat Steps 1 and 2 $B$ times and obtain after ranking $L_n^{1*} \leq L_n^{2*} \leq \cdots \leq L_n^{B*}$.

Let $l_\alpha^*$ be the $1 - \alpha$ quantile of $L_n^*$ satisfying $P(L_n^* \geq l_\alpha^* | \{X_t\}_{t=1}^{n+1}) = \alpha$. A Monte Carlo approximation of $l_\alpha^*$ is $L_n^{[B(1-\alpha)]+1*}$. The proposed test rejects $H_0$ if $L_n \geq l_\alpha^*$.

The following theorem is the bootstrap version of Theorem 1 establishing the convergence in joint distribution of the bootstrap version of the test statistics.

THEOREM 2. *Under Assumptions 1–4, as $n \to \infty$, given $\{X_t\}_{t=1}^{n+1}$, the conditional distribution of $L_n^*$ converges to the distribution of $\max_{1 \leq k \leq J} Z_k$ in probability as $n \to \infty$, where $(Z_1, \ldots, Z_J)^T \sim N(\beta \mathbf{1}_J, \Sigma_J)$.*

The next theorem shows that the proposed EL test based on the bootstrap calibration has correct size asymptotically under $H_0$ and is consistent under $H_1$.

THEOREM 3. *Under Assumptions 1–4, $\lim_{n \to \infty} P(L_n \geq l_\alpha^*) = \alpha$ under $H_0$; and $\lim_{n \to \infty} P(L_n \geq l_\alpha^*) = 1$ under $H_1$.*

**5. Computation.** The computation of the proposed EL test statistic $N(h)$ involves first computing the local EL ratio $\ell\{\tilde{p}_{\hat{\theta}}(y|x)\}$ over a grid of $(x, y)$-points within the set $S \subset \mathcal{X}^2$. The number of grid points should be large enough to ensure good approximation by the Riemann sum. On top of this is the bootstrap procedure that replicates the above computation for a large number of times. The most time assuming component of the computation is the nonlinear optimization carried out when obtaining the local EL ratio $\ell\{\tilde{p}_{\hat{\theta}}(y|x)\}$. The combination of the EL and bootstrap makes the computation intensive.

In the following, we consider using a simpler version of the EL, the least squares empirical likelihood (LSEL), to formulate the test statistic. The log LSEL ratio evaluated at $\tilde{p}_{\hat{\theta}}(y|x)$ is

$$\ell^{ls}\{\tilde{p}_{\hat{\theta}}(y|x)\} = \min \sum_{t=1}^{n} \{nq_t(x, y) - 1\}^2$$

subject to $\sum_{t=1}^{n} q_t(x, y) = 1$ and $\sum_{t=1}^{n} q_t(x, y)T_t(x, y) = 0$. Let $T(x, y) = \sum_{t=1}^{n} T_t(x, y)$ and $S(x, y) = \sum_{t=1}^{n} T_t^2(x, y)$. The LSEL is much easier to compute as there are closed-form solutions for the weights $q_t(x, y)$ and hence avoids the expensive nonlinear optimization of the EL computation. According to Brown and Chen [8], the LSEL weights $q_t(x, y) = n^{-1} + \{n^{-1}T(x, y) - T_t(x, y)\}^{\tau}S^{-1}(x)T(x, y)$ and

$$\ell^{ls}\{\tilde{m}_{\hat{\theta}}(x)\} = S^{-1}(x, y)T^2(x, y),$$

which is readily computable. The LSEL counterpart to $N(h)$ is

$$N^{ls}(h) = \iint \ell^{ls}\{\tilde{m}_{\hat{\theta}}(x, y)\}\omega(x, y)\,dx\,dy$$

and the final test statistic $L_n$ becomes $\max_{h \in \mathcal{H}}(\sqrt{2}h)^{-1}\{N^{ls}(h) - 1\}$. It can be shown from Brown and Chen [8] that $N^{ls}(h)$ and $N(h)$ are equivalent to the first order. Therefore, those first-order results in Theorems 1, 2 and 3 continue to hold for the LSEL formulation. One may just use this less expensive LSEL to carry out the testing. In fact, the least squares EL formulation was used in all the simulation studies reported in the next section.

One may use the leading order term in the expansion of the EL test statistic as given in (3.7) as the local test statistic. However, doing so would require estimation of secondary "parameters," like $p(y|x)$ and $\pi(x)$. The use of LSEL or EL avoids the secondary estimation as they Studentize automatically via their respective optimization procedures.

**6. Simulation studies.** We report results of simulation studies which were designed to evaluate the empirical performance of the proposed EL test. To gain information on its relative performance, Hong and Li's test is performed for the same simulation.

Throughout the paper, the biweight kernel $K(u) = \frac{15}{16}(1 - u^2)^2 I(|u| \leq 1)$ was used in all the kernel estimation. In the simulation, we set $\Delta = \frac{1}{12}$, implying monthly observations which coincide with that of the Federal fund rate data to be analyzed. We chose $n = 125, 250$ or $500$ respectively corresponding roughly to 10 to 40 years of data. The number of simulations was 500 and the number of bootstrap resample paths was $B = 250$.

6.1. *Size evaluation.* Two simulation studies were carried out to evaluate the size of the proposed test for both the Vasicek and Cox, Ingersoll and Ross (CIR) [13] diffusion models.

6.1.1. *Vasicek models.* We first consider testing Vasicek model

$$dX_t = \kappa(\alpha - X_t)\,dt + \sigma\,dB_t.$$

The vector of parameters $\theta = (\alpha, \kappa, \sigma^2)$ takes three sets of values which correspond to Model $-2$, Model 0 and Model 2 of Pritsker [42]. The baseline Model 0 assigns $\kappa_0 = 0.85837$, $\alpha_0 = 0.089102$ and $\sigma_0^2 = 0.0021854$ which matches estimates of Aït-Sahalia [1] for an interest rate data. Model $-2$ is obtained by quadrupling $k_0$ and $\sigma_0^2$ and Model 2 by halving $k_0$ and $\sigma_0^2$ twice while keeping $\alpha_0$ unchanged. The three models have the same marginal distribution $N(\alpha_0, V_E)$, where $V_E = \frac{\sigma^2}{2\kappa} = 0.001226$. Despite the stationary distribution being the same, the models offer different levels of dependence as quantified by the mean-reverting parameter $\kappa$. From Models $-2$ to 2, the process becomes more dependent as $\kappa$ gets smaller.

The region $S$ was chosen based on the underlying transitional density so that the region attained more than 90% of the probability. This is consistent with our earlier recommendation to choose $S$ based on the kernel estimate of the transitional density. In particular, for Models $-2$, 0 and 2, it was chosen by rotating respectively $[0.035, 0.25] \times [-0.03, 0.03]$, $[0.03, 0.22] \times [-0.02, 0.02]$ and $[0.02, 0.22] \times [-0.009, 0.009]$ 45 degrees anti clock-wise. The weight function $\omega(x, y) = |S|^{-1}I\{(x, y) \in S\}$, where $|S|$ is the area of $S$.

Both the cross-validation (CV) and the reference to a bivariate normal distribution (the Scott Rule, Scott [45]) method were used to select the bandwidth set $\mathcal{H}$. A table in a full report to this paper (Chen, Gao and Tang [11]) reports the average bandwidths obtained by the two methods. We observed that, for each given $n$, regardless of which method was used, the chosen bandwidth became smaller as the model was shifted from Model $-2$ to Model 2. This indicated that both methods took into account the changing level of dependence induced by these models. We considered two methods in choosing the bandwidth set. One was to choose six bandwidths for each combination of model and sample size that contained the average $h_{cv}$ within the lower range of $\mathcal{H}$. The second approach was to select $h_{ref}$ given the Scott Rule for each sample and then choose other $h$-values in the set by setting $a = 0.95$ so that $h_{ref}$ is the third smallest bandwidth in the set of six. This second approach could be regarded as data-driven as it was different from sample to sample.

The maximum likelihood estimator was used to estimate $\theta$ in each simulation and each resample in the bootstrap. Again, a table in Chen, Gao and Tang [11] summarizes the quality of the parameter estimation, which showed that the estimation of $\kappa$ was subject to severe bias when the mean-reversion is weak. The deterioration in the quality of the estimates, especially for $\kappa$ when the dependence became stronger, was quite alarming.

The average sizes of the proposed test at the nominal size of 5% using the two bandwidth set selection rules are reported in Table 1. It shows that the sizes of the proposed test using the proposed two bandwidth set selections were quite close to the nominal level consistently for the sample sizes considered. For Model 2, which has the weakest mean-reversion, there was some size distortion when $n = 125$ for the fix bandwidth selection. However, it was significantly alleviated when $n$ was increased. The message conveyed by Table 1 is that we need not have a large number of years of data in order to achieve a reasonable size for the test. Table 1 also reports the single-bandwidth based test based on $N(h)$ and the asymptotic normality as conveyed by Theorem 1 with $J = 1$. However, the asymptotic test has severe size distortion and highlights the need for the bootstrap procedure.

We then carried out simulation for the test of Hong and Li [32]. The Scott Rule adopted by Hong and Li was used to get an initial bandwidth $h_{scott} = \hat{S}_z n^{-1/6}$, where $\hat{S}_z$ is the sample standard deviation of the transformed series. There was little difference in the average of $h_{scott}$ among the three Vasicek models in the simulation. We used the one corresponding to Vasicek $-2$. We then chose 2 equally spaced bandwidths below and above the average $h_{scott}$. The nominal 5% test at each bandwidth was carried out with the lag value 1. For the sample sizes considered, the sizes of the test did not settle well at the nominal level, similar to what happened for the asymptotic test as reported in Table 1. We then carried out the proposed parametric bootstrap procedure for Hong and Li's test. As shown in Table 2, the bootstrap largely improved the size of the test.

### 6.1.2. *CIR models.*
We then conduct simulation on the CIR process

$$dX_t = \kappa(\alpha - X_t)\,dt + \sigma\sqrt{X_t}\,dB_t \tag{6.1}$$

to see if the pattern of results observed for the Vasicek models holds for the CIR models. The parameters were the following: $\kappa = 0.89218$, $\alpha = 0.09045$ and $\sigma^2 = 0.032742$ in the first model (CIR 0); $\kappa = 0.44609$, $\alpha = 0.09045$

Table 1

*Empirical sizes (in percentage) of the proposed EL test (the last two columns) and the single bandwidth based test (in the middle) for the Vasicek models, as well as those of the single bandwidth test based on the asymptotic normality (in round bracket): $\alpha_1$—size based on the fixed bandwidth set; $\alpha_2$—size based one the data-driven bandwidths*

| | **Bandwidths** | | | | | | $\alpha_1$ | $\alpha_2$ |
|---|---|---|---|---|---|---|---|---|
| | | | A: Model $-2$ | | | | | |
| $n = 125$ | 0.030 | 0.032 | 0.034 | 0.036 | 0.0386 | 0.041 | | |
| Size | 9.4 | 8.2 | 5.2 | 4.6 | 3 | 2.4 | 4.4 | 7.6 |
| | (40.4) | (38.8) | (34.8) | (34.8) | (34.2) | (34.8) | | |
| $n = 250$ | 0.022 | 0.023 | 0.024 | 0.026 | 0.0269 | 0.0284 | | |
| Size | 8 | 5.2 | 4.6 | 4.4 | 3.4 | 2.6 | 4.6 | 6 |
| | (34.4) | (28.6) | (24) | (21) | (17.4) | (16) | | |
| $n = 500$ | 0.02 | 0.021 | 0.022 | 0.023 | 0.0245 | 0.0258 | | |
| Size | 6.2 | 5.8 | 5.4 | 5.2 | 5 | 5 | 5.4 | 5.6 |
| | (29.6) | (23.6) | (19) | (14.6) | (10.8) | (8.4) | | |
| | | | B: Model 0 | | | | | |
| $n = 125$ | 0.016 | 0.017 | 0.019 | 0.020 | 0.022 | 0.024 | | |
| Size | 5.8 | 6 | 6 | 4.2 | 4.4 | 3 | 4.2 | 7 |
| | (43) | (39) | (36.6) | (34.8) | (36.6) | (37) | | |
| $n = 250$ | 0.014 | 0.015 | 0.017 | 0.018 | 0.02 | 0.022 | | |
| Size | 6 | 6.2 | 6.2 | 3.8 | 2.4 | 2.8 | 5.2 | 5.8 |
| | (31.6) | (27) | (20.6) | (20) | (17.8) | (17.8) | | |
| $n = 500$ | 0.01 | 0.011 | 0.012 | 0.013 | 0.015 | 0.016 | | |
| Size | 6.8 | 4.4 | 5.2 | 6.4 | 5.6 | 4 | 5.4 | 6.2 |
| | (36.4) | (26.8) | (20.6) | (13) | (11.2) | (9.2) | | |
| | | | C: Model 2 | | | | | |
| $n = 125$ | 0.008 | 0.009 | 0.010 | 0.011 | 0.013 | 0.014 | | |
| Size | 12.6 | 11 | 10 | 14.6 | 14.4 | 13.6 | 12.6 | 3.4 |
| | (60) | (53.4) | (47.2) | (46) | (45) | (42.2) | | |
| $n = 250$ | 0.006 | 0.007 | 0.008 | 0.009 | 0.01 | 0.011 | | |
| Size | 12.2 | 10 | 7.4 | 8.8 | 7 | 11 | 8.8 | 4.2 |
| | (39) | (35) | (31) | (30) | (31) | (33.2) | | |
| $n = 500$ | 0.004 | 0.005 | 0.0054 | 0.0063 | 0.0074 | 0.0086 | | |
| Size | 8.2 | 8.4 | 8 | 8.6 | 7 | 9 | 7.2 | 5.6 |
| | (75.2) | (63) | (51.6) | (39.8) | (32.8) | (24.4) | | |

and $\sigma^2 = 0.016371$ in model CIR 1 and $\kappa = 0.22305$, $\alpha = 0.09045$ and $\sigma^2 = 0.008186$ in model CIR 2. CIR 0 was the model used in Pritsker [42] for power evaluation. The region $S$ was chosen by rotating 45-degrees anti-clockwise $[0.015, 0.25] \times [-0.015, 0.015]$ for CIR 0, $[0.015, 0.25] \times [-0.012, 0.012]$ for CIR 1 and $[0.015, 0.25] \times [-0.008, 0.008]$ for CIR 2, respectively. All the regions have a coverage probability of at least 0.90.

Table 3 reports the sizes of the proposed test based on two bandwidth sets, as well as the single bandwidth-based tests that involve the bootstrap

TABLE 2
*Asymptotic and bootstrap (in parentheses) sizes of Hong and Li's test for the Vasicek models*

| $n = 125$ | $h$: | 0.08 | 0.1 | 0.12 | 0.14 | 0.16 |
|---|---|---|---|---|---|---|
| Vasicek $-2$ | | 11.8 (6.2) | 3.8 (5.8) | 1.2 (5.2) | 0.6 (4.8) | 0.6 (4.4) |
| Vasicek 0 | | 13.2 (4.6) | 4.2 (4.2) | 1.2 (3.8) | 0.8 (3.4) | 0.8 (3.4) |
| Vasicek 2 | | 11.8 (4.4) | 4.2 (3.2) | 2 (3.2) | 1.4 (2.8) | 1.6 (2.8) |
| $n = 250$ | $h$: | 0.07 | 0.09 | 0.11 | 0.13 | 0.15 |
| Vasicek $-2$ | | 15 (6.2) | 6.2 (5.4) | 2.2 (6) | 1.6 (5.2) | 1 (6.4) |
| Vasicek 0 | | 13 (4.8) | 5.8 (5.4) | 3.4 (5.4) | 1.8 (5.8) | 2 (6.6) |
| Vasicek 2 | | 15 (5.4) | 7.4 (5.4) | 3 (5.8) | 1.6 (6.6) | 1.4 (5.4) |
| $n = 500$ | $h$: | 0.06 | 0.08 | 0.1 | 0.12 | 0.14 |
| Vasicek $-2$ | | 15.6 (5.6) | 7.6 (5.6) | 2.8 (4.6) | 2.4 (4.8) | 1.2 (6.2) |
| Vasicek 0 | | 19.4 (5.4) | 8.4 (5.6) | 3.4 (5.4) | 2.6 (6.2) | 1.6 (6.6) |
| Vasicek 2 | | 17 (6.6) | 9 (5.8) | 4.6 (5.4) | 3.2 (6.6) | 2.4 (7.6) |

TABLE 3
*Empirical sizes (in percentage) of the proposed EL test (the last two columns) and the single bandwidth based test (in the middle) for the CIR models: $\alpha_1$—size based on the fixed bandwidth set; $\alpha_2$—size based one the data-driven bandwidths*

| | **Bandwidths** | | | | | | $\alpha_1$ | $\alpha_2$ |
|---|---|---|---|---|---|---|---|---|
| | A: CIR 0 | | | | | | | |
| $n = 125$ | 0.022 | 0.025 | 0.029 | 0.033 | 0.038 | 0.044 | | |
| Size | 4.8 | 4.8 | 4.8 | 3.2 | 2.4 | 2 | 3.0 | 6.6 |
| $n = 250$ | 0.018 | 0.021 | 0.024 | 0.028 | 0.032 | 0.037 | | |
| Size | 5.2 | 5.6 | 5.2 | 5.2 | 4 | 3.8 | 5.0 | 6 |
| $n = 500$ | 0.016 | 0.018 | 0.021 | 0.024 | 0.027 | 0.031 | | |
| Size | 4.2 | 5.4 | 4.6 | 4.8 | 3.6 | 4.2 | 4.8 | 5.6 |
| | B: CIR 1 | | | | | | | |
| $n = 125$ | 0.017 | 0.02 | 0.022 | 0.026 | 0.03 | 0.035 | | |
| Size | 5.4 | 3.6 | 4.0 | 4.2 | 3.2 | 2.6 | 3.8 | 7.6 |
| $n = 250$ | 0.014 | 0.016 | 0.018 | 0.021 | 0.024 | 0.028 | | |
| Size | 5.2 | 6.6 | 5.6 | 5.6 | 5.6 | 3.4 | 5.2 | 5.6 |
| $n = 500$ | 0.012 | 0.014 | 0.016 | 0.018 | 0.021 | 0.024 | | |
| Size | 5.4 | 3.8 | 4.4 | 4.4 | 4.6 | 4 | 5.2 | 5.2 |
| | C: CIR 2 | | | | | | | |
| $n = 125$ | 0.012 | 0.014 | 0.016 | 0.018 | 0.021 | 0.024 | | |
| Size | 7.6 | 7.2 | 7.6 | 6.2 | 6.6 | 4.2 | 6.8 | 8.4 |
| $n = 250$ | 0.01 | 0.012 | 0.013 | 0.015 | 0.017 | 0.02 | | |
| Size | 5.6 | 6.4 | 5.8 | 6.8 | 6.2 | 5.4 | 6.2 | 7.6 |
| $n = 500$ | 0.008 | 0.009 | 0.011 | 0.012 | 0.014 | 0.016 | | |
| Size | 4 | 4.2 | 3.6 | 4 | 4.8 | 3.4 | 4 | 6.6 |

Table 4A

*Empirical power (in percentage) of the proposed EL test (last two columns) and the single bandwidth based test: $\alpha_1$—power based on the fixed bandwidth set; $\alpha_2$—power based on the data-driven bandwidths*

| $n$ | | Single bandwidth-based tests | | | | | $\alpha_1$ | $\alpha_2$ |
|-----|-----|--------|--------|--------|--------|--------|-------|-------|
| 125 | $h$ | 0.0199 | 0.0219 | 0.0241 | 0.0265 | 0.0291 | | |
| | Power | 80.4 | 74 | 67.2 | 66.4 | 65.2 | 79.8 | 63.6 |
| 250 | $h$ | 0.0141 | 0.0158 | 0.0177 | 0.0199 | 0.0223 | | |
| | Power | 87.6 | 81.2 | 76.4 | 74 | 72.8 | 88.6 | 65.2 |
| 500 | $h$ | 0.0113 | 0.0126 | 0.0141 | 0.0157 | 0.0175 | | |
| | Power | 90.8 | 84.8 | 82.8 | 84.4 | 80.8 | 96.8 | 81.4 |

simulation. The bandwidth sets were chosen based on the same principle as outlined for the Vasicek models and are reported in the table. We find that the proposed test continued to have reasonable size for the three CIR models despite that there were severe biases in the estimation of $\kappa$. The size of the single bandwidth based tests as well as the overall test were quite respectable for $n = 125$. It is interesting to see that despite $\kappa$ still being poorly estimated for the CIR 2, the severe size distortion observed earlier for Vasicek 2 for the fixed bandwidth set was not present. Hong and Li's test was also performed for the three CIR models. The performance was similar to that of the Vasicek models reported in Table 2 and, hence, we would not report here.

6.2. *Power evaluation.* To gain information on the power of the proposed test, we carried out simulation to test for the Vasicek model, while the real process was the CIR 0 as in Pritsker's power evaluation of Aït-Sahalia's test. The region $S$ was obtained by rotating $[0.015, 0.25] \times [-0.015, 0.015]$ 45 degrees anti-clock wise. The average CV bandwidths based on 500 simulations were 0.0202 (the standard error of 0.0045) for $n = 125$, 0.016991 (0.00278) for $n = 250$ and 0.014651 (0.00203) for $n = 500$.

Table 4B

*Asymptotic and bootstrap power (in round brackets) of Hong and Li's test*

| $n$ | | Asymptotic power (bootstrap power) | | | | |
|-----|-----|-----------|----------|-----------|-----------|-----------|
| 125 | $h$ | 0.08 | 0.1 | 0.12 | 0.14 | 0.16 |
| | Power | 26 (4.6) | 16.4 (3.8) | 9.6 (2.8) | 5.8 (3) | 4.8 (3.6) |
| 250 | $h$ | 0.07 | 0.09 | 0.11 | 0.13 | 0.15 |
| | Power | 41 (5.8) | 29 (6) | 18.8 (5.4) | 14.6 (7.2) | 10.4 (6.4) |
| 500 | $h$ | 0.06 | 0.08 | 0.1 | 0.12 | 0.14 |
| | Power | 57.4 (6) | 49 (5.4) | 40.2 (5) | 34 (6.2) | 31.2 (7.2) |

Table [4A] reports the power of the EL test and the single bandwidth-based tests, including the fixed bandwidth sets used in the simulation. We find the tests had quite good power. As expected, the power increased as $n$ increased. One striking feature was that the power of the test tends to be larger than the maximum power of the single bandwidth-based tests, which indicates that it is worthwhile to formulate the test based on a set of bandwidths. Table [4B] also reports the power of Hong and Li's test. It is found that while the bootstrap calibration improved the size of the test, it largely reduced the power. The power reduction was quite alarming. In some cases, the test had little power. Our simulation results on Hong and Li's test were similar to those reported in Aït-Sahalia, Fan and Peng [5].

**7. Case studies.** We apply the proposed test on the Federal fund rate data set between January 1963 and December 1998 which has $n = 432$ observations. Aït-Sahalia [2] used this data set to demonstrate the performance of the maximum likelihood estimation. We test for five popular one-factor diffusion models which have been proposed to model interest rate dynamics. In additional to the Vasicek and CIR processes, we consider

$$(7.1) \qquad dX_t = X_t\{\kappa - (\sigma^2 - \kappa\alpha)X_t\}\,dt + \sigma X_t^{3/2}\,dB_t,$$

$$(7.2) \qquad dX_t = \kappa(\alpha - X_t)\,dt + \sigma X_t^\rho\,dB_t,$$

$$(7.3) \qquad dX_t = (\alpha_{-1}X_t^{-1} + \alpha_0 + \alpha_1 X_t + \alpha_2 X_t^2)\,dt + \sigma X_t^{3/2}\,dB_t.$$

They are respectively the inverse of the CIR process (ICIR) (7.1), the constant elasticity of the volatility (CEV) model (7.2) and the nonlinear drift (NL) model (7.3) of Aït-Sahalia [1].

The data are displayed in Figure [1](a), which indicates a strong dependence as they scattered around a narrow band around the 45-degree line. There was an increased volatility when the rate was larger than 12%. The model-implied transitional densities under the above five diffusion models are displayed in the other panels of Figure [1] using the MLEs given in Aït-Sahalia [2], which were also used in the formulation of the proposed test statistic. Figure [1] shows that the densities implied by the Inverse CIR, the CEV and the nonlinear drift models were similar to each other, and were quite different from those of the Vasicek and CIR models. The bandwidths prescribed by the Scott rule and the CV for the kernel estimation were respectively $h_{\mathrm{ref}} = 0.007616$ and $h_{\mathrm{cv}} = 0.00129$. Plotting the density surfaces indicated that a reasonable range for $h$ was from 0.007 to 0.02, which offered a lot of smoothness from slightly undersmoothing to slightly oversmoothing. This led to a bandwidth set consisting of $J = 7$ bandwidths with $h_1 = 0.007$, $h_J = 0.020$ and $a = 0.8434$.

Kernel transitional density estimates and the smoothed model-implied transitional densities for the five models are plotted in Figure [2] for $h =$ 0.007. By comparing Figure 2 with Figure 1, we notice the effect of kernel smoothing on these model-implied densities. In formulating the final test statistic $L_n$, we chose

$$(7.4) \qquad N(h) = \frac{1}{n} \sum_{t=1}^{n} \ell \{ \tilde{p}_{\hat{\theta}}(X_{t+1} | X_t) \} \omega_1(X_t, X_{t+1}),$$

where $\omega_1$ is a uniform weight over a region by rotating $[0.005, 0.4] \times [-0.03, 0.03]$ 45 degrees anti clock-wise. The region contains all the data values of the pair $(X_t, X_{t+1})$. As seen from (7.4), $N(h)$ is asymptotically equivalent to the statistic defined in (3.8) with $\omega(x, y) = p(x, y)\omega_1(x, y)$.

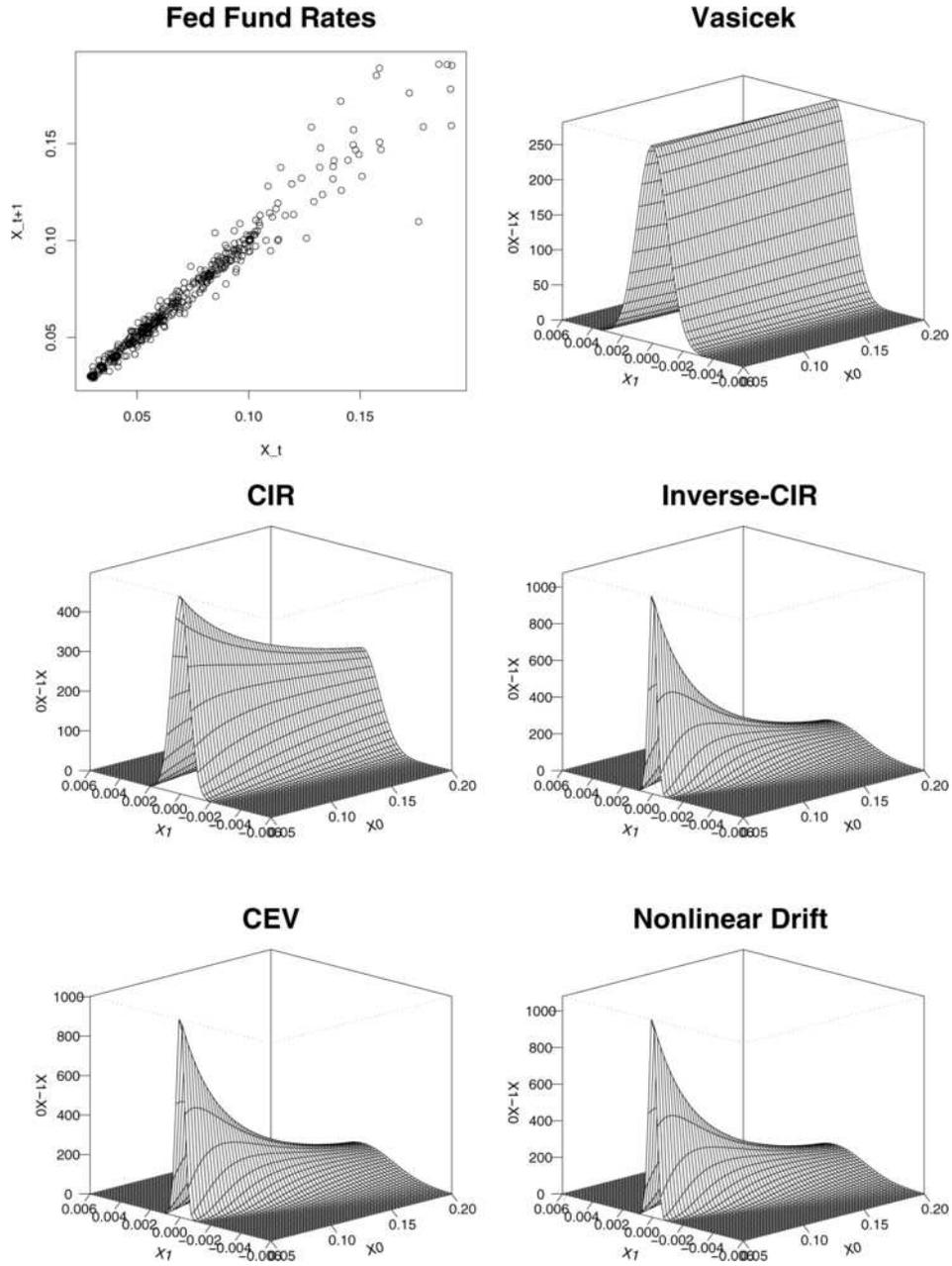

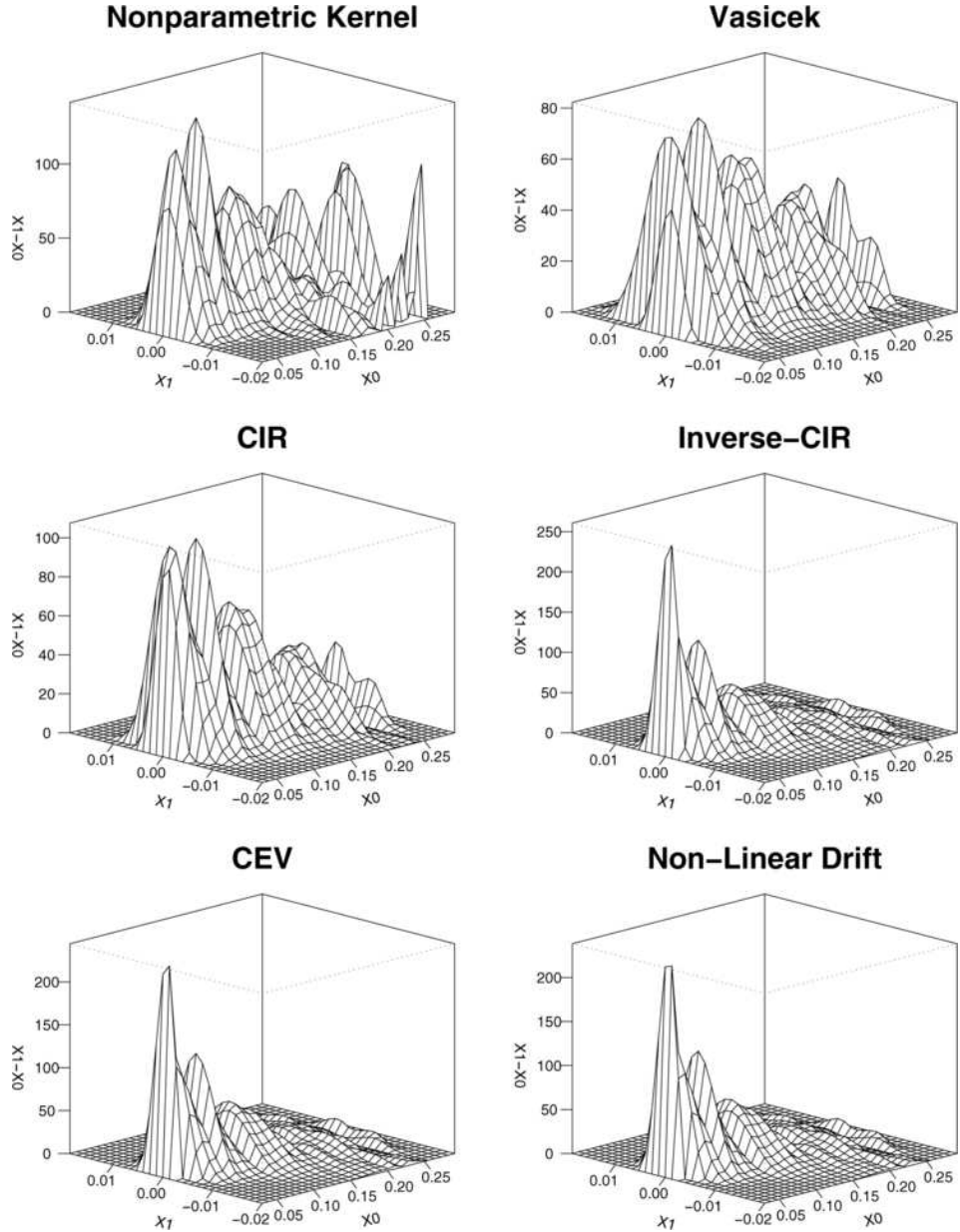

FIG. 2.  *Nonparametric kernel transitional density estimate and smoothed parametric transitional densities.*

The $p$-values of the proposed tests are reported in Table 5, which were obtained based on 500 bootstrap resamples. It shows little empirical support for the Vasicek model and quite weak support for the CIR. What was surprising is that there was some empirical support for the inverse CIR, the CEV and the nonlinear drift models. In particular, for CEV and the nonlinear drift models, the $p$-values of the single bandwidth based tests were all quite supportive even for small bandwidths. Indeed, by looking at Figure 2, we see quite noticeable agreement between the nonparametric kernel density estimates and the smoothed densities implied by the CEV and nonlinear drift models.

**8. Conclusion.** The proposed test shares some similar features with the test proposed in Hong and Li [32]. For instance, both are applicable to test continuous-time and discrete-time Markov processes by focusing on the specification of the transitional density. An advantage of Hong and Li's test is its better handling of nonstationary processes. The proposed test is based on a direct comparison between the kernel estimate and the smoothed model-implied transitional density, whereas Hong and Li's test is an indirect comparison after the probability integral transformation. An advantage of the direct approach is its robustness against poor quality parameter estimation which is often the case for weak mean-reverting diffusion models. This is because both the shape and the orientation of the transitional density are much less affected by the poor quality parameter estimation. Another aspect is that Hong and Li's test is based on asymptotic normality and can be under the influence of slow convergence despite the fact that the transformed series is asymptotically independent. Indeed, our simulation showed that it is necessary to implement the bootstrap procedure for Hong and Li's test. The last and the most important aspect is that a test based on the conditional distribution transformation tends to reduce the power comparing with a direct test based on the transitional density. This has been indicated by our simulation study, as well as that of Aït-Sahalia, Fan and Peng [5]. Interested readers can read Aït-Sahalia, Fan and Peng [5] for more insights.

Our proposed test is formulated for the univariate diffusion process. An extension to multivariate diffusion processes can be made by replacing the univariate kernel smoothing in estimating the transitional density with multivariate kernel transitional density estimation. There is no substantial difference in the formulation of the EL test statistic and the parametric bootstrap procedure. If the exact form of the transitional density is unknown, which is more likely for multivariate diffusion processes, the approximate transitional density expansion of Aït-Sahalia and Kimmel [6] is needed.

Table 5
*P-values for the federal fund rate data*

| Model | Vasicek | CIR | ICIR | CEV | NL |
|---|---|---|---|---|---|
| Test statistics $L_n$ | 29.71 | 12.80 | 66.63 | 64.56 | 69.10 |
| Critical value $l^*_{0.05}$ | 2.54 | 22.27 | 303.4 | 434.77 | 557.52 |
| $p$-value | 0.0 | 0.142 | 0.294 | 0.434 | 0.422 |

## APPENDIX

As the Lagrange multiplier $\lambda(x,y)$ is implicitly dependent on $h$, we need first to extend the convergence rate for a single $h$-based $\sup_{(x,y) \in S} \lambda(x,y)$ conveyed in (3.6) to be valid uniformly over $\mathcal{H}$. To prove Theorem 1, we need the following lemmas first.

LEMMA A.1. *Under Assumptions 1–4,*

$$\max_{h \in \mathcal{H}} \sup_{(x,y) \in S} \lambda(x,y) = o_p\{n^{-1/4}\log(n)\}.$$

PROOF. For any $\delta > 0$,

$$P\left(\max_{h \in \mathcal{H}_n} \sup_{(x,y) \in S} h\lambda(x,y) \geq \delta n^{-1/2}\log(n)\right)$$

$$\leq \sum_{h \in \mathcal{H}} P\left(\sup_{(x,y) \in S} h\lambda(x,y) \geq \delta n^{-1/2}\log(n)\right).$$

As the number of bandwidths in $H$ is finite, by checking the relevant derivations in Chen, Härdle and Li [12], it can be shown that

$$P\left(\sum_{(x,y) \in S} h\lambda(x,y) \geq \delta n^{-1/2}\log(n)\right) \to 0$$

as $n \to \infty$. This implies that $\max_{h \in \mathcal{H}} \sup_{(x,y) \in S} h\lambda(x,y) = o_p\{\delta n^{-1/2}\log(n)\}$. Then the lemma is established by noting that $h_1$, the smallest bandwidth in $\mathcal{H}$, is of order $n^{-\gamma_1}$, where $\gamma_1 \in (1/7, 1/4)$ as assumed in Assumption 1. □

Before introducing more lemmas, we present some expansions for the EL test statistic $N(h)$. Let

$$\tilde{p}_\theta(x,y) = \tilde{p}_\theta(y|x)\hat{\pi}(x)$$

and

$$\tilde{p}(x,y) = n^{-1} \sum_{t=1}^{n+1} K_h(x - X_t) \sum_{s=1}^{n+1} w_s(y) p(X_s|X_t)$$

be the kernel smoothed versions of the parametric and nonparametric joint densities $p_\theta(x,y)$ and $p(x,y)$, respectively. Due to the relationship between transitional and joint densities,

$$
\begin{aligned}
N(h) &= (nh^2) \iint \frac{\{\hat{p}(x,y) - \tilde{p}_{\hat{\theta}}(x,y)\}^2}{R^2(K)p(x,y)} \omega(x,y)\, dx\, dy \\
&\quad + \tilde{O}_p\{h^2 + (nh^2)^{-1/2}\log^3(n)\} \\
&= (nh^2)R^{-2}(K) \\
&\quad \times \iiint \Bigg[ \frac{\{\hat{p}(x,y) - \tilde{p}(x,y)\}^2}{p(x,y)} \\
&\qquad + \frac{2\{\hat{p}(x,y) - \tilde{p}(x,y)\}\{\tilde{p}(x,y) - \tilde{p}_{\hat{\theta}}(x,y)\}}{p(x,y)} \\
&\qquad + \frac{\{\tilde{p}(x,y) - \tilde{p}_{\hat{\theta}}(x,y)\}^2}{p(x,y)} \Bigg] \\
&\quad \times \omega(x,y)\, dx\, dy + \tilde{O}_p\{h^2 + (nh^2)^{-1/2}\log^3(n)\} \\
&=: N_1(h) + N_{2\hat{\theta}}(h) + N_{3\hat{\theta}}(h) + \tilde{O}_p\{h^2 + (nh^2)^{-1/2}\log^3(n)\}.
\end{aligned}
$$

(A.1)

Here and throughout the proofs, $\tilde{o}(\delta_n)$ and $\tilde{O}(\delta_n)$ denote stochastic quantities which are respectively $o(\delta_n)$ and $O(\delta_n)$ uniformly over $S$ for a nonnegative sequence $\{\delta_n\}$.

Using Assumptions 3 and 4, we have $N_{l\hat{\theta}}(h) = N_{l\theta^*}(h) + \tilde{o}_p(h)$, where $\theta^* = \theta_0$ under $H_0$ and $\theta_1$ under $H_1$. Thus,

$$
\begin{aligned}
N(h) &= N_1(h) + N_{2\theta^*}(h) + N_{3\theta^*}(h) + \tilde{o}_p(h) \\
&\quad + \tilde{O}_p\{(nh^2)^{-1/2}\log^3(n)\}.
\end{aligned}
$$

We start with some lemmas on $\hat{p}(x,y)$, $\tilde{p}(x,y)$ and $\tilde{p}_\theta(x,y)$. Let $K^{(2)}$ be the convolution of $K$, $MK^{(2)}(t) = \iint uK(u)K(t+u)\, du$ and $p_3(x,y,z)$ be the joint density of $(X_t, X_{t+1}, X_{t+2})$.

The following lemmas are presented without proofs. The detailed proofs are given in Chen, Gao and Tang [11].

LEMMA A.2.  *Under Assumptions 1–4,*

$$\text{Cov}\{\hat{p}(s_1, t_1), \hat{p}(s_2, t_2)\}$$

$$= \frac{K^{(2)}((s_2 - s_1)/h)K^{(2)}((s_2 - s_1)/h)p(s_1, t_1)}{nh^2}$$

$$- \left( \frac{MK^{(2)}((s_2 - s_1)/h)\,\partial p(s_1, t_1)/\partial x}{nh} \right.$$

$$\left. + \frac{MK^{(2)}((t_2 - t_1)/nh)\,\partial p(s_1, t_1)/\partial y}{nh} \right)$$

$$+ \frac{p_3(s_1, t_1, t_2)K^{(2)}((s_2 - t_1)/h) + p_3(s_2, t_2, t_1)K^{(2)}((s_1 - t_2)/h)}{nh}$$

$$+ o\{(nh)^{-1}\}.$$

LEMMA A.3.  *Suppose that Assumptions 1–4 hold. Let* $\Delta_\theta(x, y) = \{p_\theta(y|x) - p(y|x)\}\pi(x)$. *Then*

$$(A.2) \quad \begin{aligned} E\{\tilde{p}_\theta(x, y) - \hat{p}(x, y)\} &= \Delta_\theta(x, y) + \frac{1}{2}h^2\sigma_K^2\left\{ \frac{\partial^2}{\partial x^2} + \frac{\partial^2}{\partial y^2} \right\}\Delta_\theta(x, y) \\ &\quad + \tilde{O}(h^3), \end{aligned}$$

$$(A.3) \quad \begin{aligned} E\{\tilde{p}_\theta(x, y) - \tilde{p}(x, y)\} &= \Delta_\theta(x, y) + \frac{1}{2}h^2\sigma_K^2\left\{ \frac{\partial^2}{\partial x^2} + \frac{\partial^2}{\partial y^2} \right\}\Delta_\theta(x, y) \\ &\quad + \tilde{O}(h^3), \end{aligned}$$

$$(A.4) \quad \begin{aligned} \text{Cov}\{\tilde{p}(s_1, t_1), \tilde{p}(s_2, t_2)\} &= \frac{K^{(2)}((t_2 - t_1)/h)p(s_1, t_1)p(s_2, t_1)}{nh\pi(t_2)} \\ &\quad + \tilde{o}\{(nh)^{-1}\}. \end{aligned}$$

LEMMA A.4.  *Under Assumptions 1–4, we have*

$$\begin{aligned} \text{Cov}\{\hat{p}(s_1, t_1), \tilde{p}(s_2, t_2)\} &= \frac{p(s_1, t_1)}{nh\pi(t_2)}\left[ K^{(2)}\left( \frac{t_2 - s_1}{h} \right)p(s_2, s_1) \right] \\ &\quad + \frac{p(s_1, t_1)}{nh\pi(t_2)}\left[ K^{(2)}\left( \frac{t_2 - s_1}{h} \right)p(s_2, s_1) \right]. \end{aligned}$$

LEMMA A.5.  *If $H_0$ is true, then $N_{2\theta^*}(h) = N_{3\theta^*}(h) = 0$ for all $h \in \mathcal{H}$.*

PROOF.  Under $H_0$, $p(y|x) = p_{\theta_0}(y|x)$ and $\theta^* = \theta_0$. Hence, $\tilde{p}(x, y) - \tilde{p}_{\theta^*}(x, y) = n^{-1}\sum K_h(x - X_t)\sum w_s(y)\{p(X_s|X_t) - p_{\theta_0}(X_s|X_t)\} = 0.$   □

Let us now study the leading term $N_1(h)$. From (A.1) and by hiding the variables of integrations,

$$
\begin{aligned}
N_1(h) &= \frac{(nh^2)}{R^2(K)} \iint \frac{\{\hat{p} - \tilde{p}\}^2}{R^2(K)p}\omega \\
&= \frac{(nh^2)}{R^2(K)} \iint \left[ \frac{\{\hat{p} - E\hat{p}\}^2}{p} + \frac{\{E\tilde{p} - \tilde{p}\}^2}{p} + \frac{\{E\hat{p} - E\tilde{p}\}^2}{p} \right. \\
&\qquad\qquad + \frac{2\{\hat{p} - E\hat{p}\}\{E\hat{p} - E\tilde{p}\}}{p} + \frac{2\{\hat{p} - E\hat{p}\}\{E\tilde{p} - \tilde{p}\}}{p} \\
&\qquad\qquad\qquad\qquad \left. + \frac{2\{E\hat{p} - E\tilde{p}\}\{E\tilde{p} - \tilde{p}\}}{p} \right]\omega \\
&=: \sum_{j=1}^{6} N_{1j}(h).
\end{aligned}
$$

We are to show in the following lemmas that $N_{11}(h)$ dominates $N_1(h)$ and $N_{1j}(h)$ for $j \geq 2$ are all negligible except $N_{12}(h)$, which contributes to the mean of $N_1(h)$ in the second order.

LEMMA A.6. *Under Assumptions 1–4, then uniformly with respect to* $\mathcal{H}$,

$$
h^{-1}E\{N_{11}(h) - 1\} = o(1), \tag{A.5}
$$

$$
\mathrm{Var}\{h^{-1}N_{11}(h)\} = \frac{2K^{(4)}(0)}{R^4(K)} \iint \omega^2(x,y)\, dx\, dy + o(1), \tag{A.6}
$$

$$
\mathrm{Cov}\{h_1^{-1}N_{11}(h_1), h_2^{-1}N_{11}(h_2)\} = \frac{2\nu(h_1/h_2)}{R^4(K)} \iint \omega^2(x,y)\, dx\, dy
$$
$$
+ o(1). \tag{A.7}
$$

PROOF. From Lemma A.1 and the fact that $MK^{(2)}(0) = 0$,

$$
\begin{aligned}
E\{N_{11}(h)\} &= \frac{nh^2}{R^2(K)} \iint \frac{\mathrm{Var}\{\hat{p}(x,y)\}}{p(x,y)}\omega(x,y)\, dx\, dy \\
&= \frac{1}{R^2(K)} \iint \left[ (K^{(2)}(0))^2 + 2hK^{(2)}\left(\frac{y-x}{h}\right)\frac{p_3(x,y,y)}{p(x,y)} \right] \\
&\qquad\qquad \times \omega(x,y)\, dx\, dy\{1 + o(1)\} \\
&= 1 + O(h^2),
\end{aligned} \tag{A.8}
$$

which leads to (A.5). To derive (A.6), let

$$
\hat{Z}_n(s,t) = (nh^2)^{1/2}\frac{\hat{p}(s,t) - E\hat{p}(s,t)}{R(K)p^{1/2}(x,t)}.
$$

It may be shown from the fact that $K$ is bounded and an other regularity condition assumed that $E\{|\hat{Z}_n(s_1,t_1)|^{2+\epsilon}|\hat{Z}_n(s_2,t_2)|^{2+\epsilon}\} \leq M$ for some positive $\epsilon$ and $M$. And hence, $\{\hat{Z}_n(s,t)\}_{n\geq1}$ and $\{\hat{Z}_n^2(s_1,t_1)\hat{Z}_n^2(s_2,t_2)\}_{n\geq1}$ are uniformly integrable respectively. Also,

$$(\hat{Z}_n(s_1,t_1),\hat{Z}_n(s_2,t_2))^T \xrightarrow{d} (Z(s_1,t_1),Z(s_2,t_2))^T,$$

which is a bivariate normal process with mean zero and covariance

$$\Sigma = \begin{pmatrix} 1 & g\{(s_1,t_1),(s_2,t_2)\} \\ g\{(s_1,t_1),(s_2,t_2)\} & 1 \end{pmatrix},$$

where $g\{(s_1,t_1),(s_2,t_2)\} = K^{(2)}(\frac{s_2-s_1}{h})K^{(2)}(\frac{t_2-t_1}{h})\frac{p^{1/2}(s_1,t_1)}{R(K)p^{1/2}(s_2,t_2)}$. Hence, by ignoring smaller order terms,

(A.9)
$$
\begin{aligned}
&\mathrm{Var}\{N_{11}(h)\} \\
&= \iiiint \mathrm{Cov}\{\hat{Z}_n^2(s_1,s_2),\hat{Z}_n^2(s_2,t_2)\} \\
&\qquad\qquad \times \omega(s_1,t_1)\omega(s_2,t_2)\,ds_1\,dt_1\,ds_2\,dt_2 \\
&= \iiiint \mathrm{Cov}\{Z^2(s_1,s_2),Z^2(s_2,t_2)\} \\
&\qquad\qquad \times \omega(s_1,t_1)\omega(s_2,t_2)\,ds_1\,dt_1\,ds_2\,dt_2 \\
&= 2\iiiint \mathrm{Cov}^2\{Z(s_1,s_2),Z(s_2,t_2)\} \\
&\qquad\qquad \times \omega(s_1,t_1)\omega(s_2,t_2)\,ds_1\,dt_1\,ds_2\,dt_2 \\
&= \frac{2}{R^4(K)}\iiiint \left\{K^{(2)}\left(\frac{s_2-s_1}{h}\right)K^{(2)}\left(\frac{t_2-t_1}{h}\right)\right\}^2 \frac{p(s_1,t_1)}{R^2(K)p(s_2,t_2)} \\
&\qquad\qquad \times \omega(s_1,t_1)\omega(s_2,t_2)\,ds_1\,dt_1\,ds_2\,dt_2 \\
&= \frac{2h^2 K^{(4)}(0)}{R^4(K)}\iint \omega^2(s,t)\pi^{-2}(t)\,ds\,dt.
\end{aligned}
$$

In the third equation above, we use a fact regarding the fourth product moments of normal random variables. Combining (A.8) and (A.9), (A.5) and (A.6) are derived. It can be checked that it is valid uniformly for all $h \in \mathcal{H}$.

The proof for (A.7) follows from that for (A.6).  $\square$

The proof of the following lemma is left in Chen, Gao and Tang [11].

LEMMA A.7.  *Under Assumptions 1–4, then uniformly with respect to* $h \in \mathcal{H}$,

(A.10)    $$h^{-1}N_{12}(h) = \frac{1}{R(K)}\iint \frac{p(x,y)}{\pi(y)}\omega(x,y)\,dx\,dy + o_p(1),$$

(A.11) $\qquad h^{-1}N_{1j}(h) = o_p(1) \qquad$ for $j \geq 3$.

Let $L(h) = \frac{1}{C(K)h}\{N(h) - 1\}$ and $\beta = \frac{1}{\sqrt{2}R(K)}\iint \frac{p(x,y)}{\pi(y)}\omega(x,y)\,dx\,dy$. In view of Lemmas A.5, A.6 and A.7, we have, under $H_0$, uniformly with respect to $\mathcal{H}$,

(A.12) $\qquad L(h) = \frac{1}{\sqrt{2}h}\{N_{11}(h) - 1\} + \beta + o_p(1)$.

Define $L_1(h) = \frac{1}{\sqrt{2}h}\{N_{11}(h) - 1\}$.

LEMMA A.8. *Under Assumptions 1–4 and $H_0$, as $n \to \infty$,*

$$(L_1(h_1), \ldots, L_1(h_J))^T \xrightarrow{d} N_J(\beta 1_J, \Sigma_J).$$

PROOF. According to the Cramér–Wold device, it suffices to show

(A.13) $\qquad \sum_{i=1}^{J} c_i L_1(h_i) \xrightarrow{d} N_J(c^\tau \beta 1_J, c^\tau \Sigma_J c)$

for an arbitrary vector of constants $c = (c_1, \ldots, c_J)^\tau$. Without loss of generality, we will only prove the case of $J = 2$. To apply Lemma A.1 of Gao and King [25], we introduce the following notation. For $i = 1, 2$, define $d_i = \frac{c_i}{\sqrt{2}h_i}$ and $\xi_t = (X_t, X_{t+1})$,

$$\epsilon_{ti}(x,y) = K\left(\frac{x - X_t}{h_i}\right)K\left(\frac{y - X_{t+1}}{h_i}\right) - E\left[K\left(\frac{x - X_t}{h_i}\right)K\left(\frac{y - X_{t+1}}{h_i}\right)\right],$$

$$\phi_i(\xi_s, \xi_t) = \frac{1}{nh_i^2}\iint \frac{\epsilon_{si}(x,y)\epsilon_{ti}(x,y)}{p(x,y)R^2(K)}\omega(x,y)\,dx\,dy,$$

$$\phi_{st} = \phi(\xi_s, \xi_t) = \sum_{i=1}^{2} d_i\phi_i(\xi_s, \xi_t) \quad \text{and} \quad \overline{L}_1(h_1, h_2) = \sum_{t=2}^{T}\sum_{s=1}^{t-1}\phi_{st}.$$

It is noted that for any given $s, t \geq 1$ and fixed $x$ and $y$, $E[\phi(x, \xi_t)] = E[\phi(\xi_s, y)] = 0$. It suffices to verify

(A.14) $\qquad \max\{M_n, N_n\}h_1^{-2} \to 0 \qquad$ as $n \to \infty$,

where

$$M_n = \max\{n^2 M_{n1}^{1/(1+\delta)}, n^2 M_{n51}^{1/(2(1+\delta))}, n^2 M_{n52}^{1/(2(1+\delta))}, n^2 M_{n6}^{1/2}\}$$

$$N_n = \max\{n^{3/2} M_{n21}^{1/(2(1+\delta))}, n^{3/2} M_{n22}^{1/(2(1+\delta))}, n^{3/2} M_{n3}^{1/2},$$

$$n^{3/2} M_{n4}^{1/(2(1+\delta))}, n^{3/2} M_{n7}^{1/(1+\delta)}\},$$

in which

$$M_{n1} = \max_{1 \le i < j < k \le n} \max\left\{ E|\phi_{ik}\psi_{jk}|^{1+\delta}, \int |\phi_{ik}\theta_{jk}|^{1+\delta}\, dP(\xi_i)\, dP(\xi_j, \xi_k) \right\},$$

$$M_{n21} = \max_{1 \le i < j < k \le n} \max\left\{ E|\phi_{ik}\phi_{jk}|^{2(1+\delta)}, \int |\phi_{ik}\phi_{jk}|^{2(1+\delta)}\, dP(\xi_i)\, dP(\xi_j, \xi_k) \right\},$$

$$M_{n22} = \max_{1 \le i < j < k \le n} \max\left\{ \int |\phi_{ik}\phi_{jk}|^{2(1+\delta)}\, dP(\xi_i, \xi_j)\, dP(\xi_k), \right.$$
$$\left. \int |\phi_{ik}\phi_{jk}|^{2(1+\delta)}\, dP(\xi_i)\, dP(\xi_j)\, dP(\xi_k) \right\},$$

$$M_{n3} = \max_{1 \le i < j < k \le n} E|\phi_{ik}\phi_{jk}|^2,$$

$$M_{n4} = \max_{1 < i,j,k \le 2n; i,j,k \text{ distinct}} \left\{ \max_P \int |\phi_{1i}\phi_{jk}|^{2(1+\delta)}\, dP \right\},$$

$$M_{n51} = \max_{1 \le i < j < k \le n} \max\left\{ E\left|\int \phi_{ik}\phi_{jk}\phi_{ik}\phi_{jk}\, dP(\xi_i)\right|^{2(1+\delta)} \right\},$$

$$M_{n52} = \max_{1 \le i < j < k \le n} \max\left\{ \int \left|\int \phi_{ik}\phi_{jk}\phi_{ik}\phi_{jk}\, dP(\xi_i)\right|^{2(1+\delta)}\, dP(\xi_j)\, dP(\xi_k) \right\},$$

$$M_{n6} = \max_{1 \le i < j < k \le n} E\left|\int \phi_{ik}\phi_{jk}\, dP(\xi_i)\right|^2,$$

$$M_{n7} = \max_{1 \le i < j < n} E[|\phi_{ij}|^{1+\delta}],$$

where the maximization in $M_{n4}$ is over the probability measures $P(\xi_1, \xi_i, \xi_j, \xi_k)$, $P(\xi_1)P(\xi_i, \xi_j, \xi_k)$, $P(\xi_1)P(\xi_{i_1})P(\xi_{i_2}, \xi_{i_3})$ and $P(\xi_1)P(\xi_i)P(\xi_j)P(\xi_k)$.

Without confusion, we replace $h_1$ by $h$ for simplicity. To verify the $M_n$ part of (A.14), we verify only

$$(A.15) \qquad\qquad \lim_{n \to \infty} n^2 h^{-2} M_{n1}^{1/(1+\delta)} = 0.$$

Let $q(x, y) = \omega(x, y)p^{-1}(x, y)$ and

$$\psi_{ij} = \frac{1}{nh^2} \int K((x - X_i)/h) K((y - X_{i+1})/h) K((x - X_j)/h)$$
$$\times K((y - X_{j+1})/h) q(x, y)\, dx\, dy$$

for $1 \le i < j < k \le n$. Direct calculation implies

$$\psi_{ik}\psi_{jk} = (nh^2)^{-2} \int K\left(\frac{x - X_i}{h}\right) K\left(\frac{y - X_{i+1}}{h}\right) K\left(\frac{x - X_k}{h}\right)$$
$$\times K\left(\frac{y - X_{k+1}}{h}\right) q(x, y) K\left(\frac{u - X_j}{h}\right) K\left(\frac{v - X_{j+1}}{h}\right)$$

$$\times K\left(\frac{u-X_k}{h}\right)K\left(\frac{v-X_{k+1}}{h}\right)q(u,v)\,dx\,dy\,du\,dv$$

$$= b_{ijk} + \delta_{ijk},$$

where $\delta_{ijk} = \psi_{ik}\psi_{jk} - b_{ijk}$ and

$$b_{ijk} = n^{-2}q(X_i, X_{i+1})q(X_j, X_{j+1})K^{(2)}\left(\frac{X_i - X_k}{h}\right)$$

$$\times K^{(2)}\left(\frac{X_j - X_k}{h}\right)K^{(2)}\left(\frac{X_{i+1} - X_{k+1}}{h}\right)K^{(2)}\left(\frac{X_{j+1} - X_{k+1}}{h}\right).$$

For any given $1 < \zeta < 2$ and $n$ sufficiently large, we may show that

$$
\begin{aligned}
M_{n11} &\le 2(E[|b_{ijk}|^\zeta] + E[|\delta_{ijk}|^\zeta]) \\
&= 2E[|b_{ijk}|^\zeta](1 + o(1)) \\
&= n^{-2\zeta}\iiint |q(x,y)q(u,v)|^\zeta \left|K^{(2)}\left(\frac{x-z}{h}\right)K^{(2)}\left(\frac{u-z}{h}\right)\right|^\zeta \\
&\qquad\qquad \times \left|K^{(2)}\left(\frac{y-w}{h}\right)K^{(2)}\left(\frac{v-w}{h}\right)\right|^\zeta \\
&\qquad\qquad \times p(x,y,u,v,z,w)\,dx\,dy\,du\,dv\,dz\,dw \\
&= C_1 n^{-2\zeta}h^4,
\end{aligned}
$$

(A.16)

where $p(x,y,u,v,z,w)$ denotes the joint density of $(X_i, X_{i+1}, X_j, X_{j+1}, X_k, X_{k+1})$ and $C_1$ is a constant. Thus, as $n \to \infty$,

(A.17) $\qquad n^2 h^{-2}M_{n11}^{1/(1+\delta)} = Cn^2 h^{-1}(n^{-2\zeta}h^2)^{1/\zeta} = h^{2(2-\zeta)/\zeta} \to 0.$

Hence, (A.17) shows that (A.15) holds for the first part of $M_{n1}$. The proof for the second part of $M_{n1}$ follows similarly. $\quad\square$

PROOF OF THEOREM 1. From (A.12) and Lemma A.8, we have, under $H_0$,

$$(L(h_1), \ldots, L(h_J)) \xrightarrow{d} N_J(\beta \mathbf{1}_J, \Sigma_J).$$

Let $Z = (Z_1, \ldots, Z_J)^T \overset{d}{\sim} N_J(\beta \mathbf{1}_J, \Sigma_J)$. By the mapping theorem, under $H_0$,

(A.18) $\qquad L_n = \max_{h \in \mathcal{H}} L(h) \xrightarrow{d} \max_{1 \le k \le J} Z_k.$

Hence, the theorem is established. $\quad\square$

Let $l_{0\alpha}$ be the upper-$\alpha$ quantile of $\max_{1 \le i \le J} Z_i$. As the distribution of $N_J(\beta \mathbf{1}_k, \Sigma_J)$ is free of $n$, so is that of $\max_{1 \le i \le J} Z_i$. And hence, $l_{0\alpha}$ is a fixed quantity with respect to $n$.

The following lemmas are required for the proof of Theorem 3.

LEMMA A.9.   *Under Assumptions 1–4, for constants $C$ and $\gamma \in (1/3, 1/2)$, $\ell\{\tilde{p}(y|x) + C(nh^2)^{-\gamma}\} \to \infty$ in probability uniformly for $(x, y) \in S$.*

PROOF.   Let $\tilde{\mu}(x, y) = \tilde{p}(y|x) + C(nh^2)^{-\gamma}$ and $Q_{t,h}(x, y) = w_{nw}, t(x) \times K_h(y - X_{t+1})$. Recall from the EL formulation given in Section 3 that $\ell\{\tilde{\mu}(x, y)\} = 2 \sum_{t=1}^{n} \log[1 + \lambda(x, y)\{Q_{t,h}(x, y) - \tilde{\mu}(x, y)\}]$, where, according to (3.4), $\lambda(x, y)$ satisfies

$$0 = \sum \frac{Q_{t,h}(x, y) - \tilde{\mu}(x, y)}{1 + \lambda(x, y)\{Q_{t,h}(x, y) - \tilde{\mu}(x, y)\}}.$$

Note that

$$
\begin{aligned}
\text{(A.19)} \quad & \lambda(x, y) \sum_{t=1}^{n} \frac{\{Q_{t,h}(x, y) - \tilde{\mu}(x, y)\}^2}{1 + \lambda(x, y)\{Q_{t,h}(x, y) - \tilde{\mu}(x, y)\}} \\
& = \sum_{t=1}^{n} \{Q_{t,h}(x, y) - \tilde{\mu}(x, y)\}.
\end{aligned}
$$

Let $S_h^2(x, y) = n^{-1} \sum_{t=1}^{n} \{Q_{t,h}(x, y) - \tilde{\mu}(x, y)\}^2$. From established results on the kernel estimator for $\alpha$-mixing sequences,

$$\text{(A.20)} \quad S_h^2(x, y) = h^{-2} R^2(K) p(y|x) p^{-1}(x) + O_p\{(nh^2)^{-2\gamma}\} + o_p(h^{-2}).$$

Note that for a positive constant $M_0$ and sufficiently large $n$, $\sup_{(x,y) \in S} |Q_{t,h}(x, y) - \tilde{\mu}(x, y)| \le h^{-2} M_0$ with probability one. Hence, (A.19) implies that

$$|\lambda(x, y) S_h^2(x, y)| \le \{1 + |\lambda(x, y)| h^{-2} M_0\} \left| n^{-1} \sum_{t=1}^{n} \{Q_{t,h}(x, y) - \tilde{\mu}(x, y)\} \right|.$$

This, along with (A.20) and the facts that

$$n^{-1} \sum_{t=1}^{n} \{Q_{t,h}(x, y) - \tilde{p}(x, y)\} = \tilde{O}_p\{(nh^2)^{-1/2} \log(n)\}$$

and $\tilde{\mu}(x, y) - \tilde{p}(y|x) = C(nh^2)^{-\gamma}$, implies

$$\text{(A.21)} \quad \lambda(x, y) = \tilde{O}_p\{h^2(nh^2)^{-\gamma}\}$$

uniformly with respect to $(x, y) \in S$. The rate of $\lambda(x, y)$ established in (A.21) allows us to carry out the Taylor expansion in (A.19) and obtain

$$\text{(A.22)} \quad \lambda(x, y) = S_h^{-2}(x, y) n^{-1} \sum_{t=1}^{n} \{Q_{t,h}(x, y) - \tilde{\mu}(x, y)\} + \tilde{O}_p\{h^4(nh^2)^{-2\gamma}\}.$$

At the same time, as $X_t$'s are continuous random variables, by applying the blocking technique and the Davydov inequality, it can be shown that, for any $\eta > 0$ and $(x, y) \in S$,

$$P(\min |Q_{t,h}(x, y) - \tilde{\mu}(x, y)| > \eta) \to 0,$$

which, using (A.19), implies

$$\text{(A.23)} \qquad |\lambda(x,y)| \geq Ch^2(nh^2)^{-\gamma}\{1 + \tilde{o}_p(1)\}.$$

From (A.22) and (A.23),

$$
\begin{aligned}
\text{(A.24)} \quad \ell\{\tilde{\mu}(x,y)\} &= 2\sum_{t=1}^{n}\log[1 + \lambda(x,y)\{Q_{t,h}(x,y) - \tilde{\mu}(x,y)\}] \\
&= n\lambda(x,y)n^{-1}\sum_{t=1}^{n}\{Q_{t,h}(x,y) - \tilde{\mu}(x,y)\} \\
&\quad + \tilde{O}_p\{nh^6(nh^2)^{-3\gamma}\} \\
&= n\lambda(x,y)\{C(nh^2)^{-\gamma} + \tilde{O}_p\{(nh^2)^{-1/2}\log n\}\} \\
&\quad + \tilde{O}_p\{h^4(nh^2)^{1-3\gamma}\} \\
&= \tilde{O}_p\{(nh^2)^{1-2\gamma}\}.
\end{aligned}
$$

As $\gamma \in (1/3, 1/2)$, $\ell\{\tilde{\mu}(x,y)\} \to \infty$ in probability as $n \to \infty$. Since both (A.21) and (A.24) are true uniformly with respect to $(x,y) \in S$, the divergence of $\ell\{\tilde{\mu}(x,y)\}$ is uniform with respect to $(x,y) \in S$. □

LEMMA A.10. *Under Assumptions 1–4 and $H_1$, for any fixed real value $x$, as $n \to \infty$, $P(L_n \geq x) \to 1$.*

PROOF. From the established theory on EL, $\ell\{\mu(x,y)\}$ is a convex function of $\mu(x,y)$, the candidate for $p(y|x)$. From the mean-value theorem and the facts that both $\ell(\cdot)$ and $\omega$ are nonnegative and $\ell$ is continuous in $(x,y)$,

$$N(h) \geq C_\delta \ell\{\tilde{p}_{\hat{\theta}}(y_0|x_0)\}$$

for some $(x_0, y_0) \in S$ and $C_\delta > 0$.

Let $\mu_1(x_0, y_0) = \tilde{p}_{\hat{\theta}}(y_0|x_0)$. By choosing $C$ properly and for $n$ large enough, we make sure that $\mu_2(x_0, y_0) = \tilde{p}(y_0|x_0) + C(nh^2)^{-\gamma}$ falls within the interval of either $(\mu_1(y_0|x_0), \hat{p}(y_0|x_0))$ or $(\hat{p}(y_0|x_0), \mu_1(y_0|x_0))$. Hence, there exists an $\alpha \in (0,1)$ such that $\mu_2(x_0, y_0) = \alpha\hat{p}(y_0|x_0) + (1-\alpha)\mu_1(y_0|x_0)$. The convexity of $\ell(\cdot)$ and the fact that $\ell\{\hat{p}(y_0|x_0)\} = 0$ lead to

$$\ell\{\mu_2(x_0, y_0)\} \leq (1-\alpha)\ell\{\mu_1(x_0, y_0)\}.$$

Since Lemma A.9 implies that $\ell\{\mu_2(x_0, y_0)\} \to \infty$ holds in probability with $\mu_2(x_0, y_0) = \mu(x_0, y_0)$, we have, as $n \to \infty$, $\ell\{\tilde{p}_{\hat{\theta}}(y_0|x_0)\} \to \infty$ in probability, which implies $N(h) \to \infty$ in probability as $n \to \infty$. This means $L(h) \to \infty$. The lemma is proved by noting $P(L_n \geq x) \geq P\{L(h_i) > x\}$ for any $i \in \{1, \ldots, J\}$. □

We now turn to the bootstrap EL test statistic $N^*(h)$, which is a version of $N(h)$ based on $\{X_t^*\}_{t=1}^{n+1}$ generated according to the parametric transitional density $p_{\tilde\theta}$. Let $\hat{p}^*(x,y)$ and $\tilde{p}_\theta^*(x,y)$ be the bootstrap versions of $\hat{p}(x,y)$ and $\tilde{p}(x,y)$ respectively, and $\tilde\theta^*$ be the maximum likelihood estimate based on the bootstrap sample. Then, the following expansion similar to (A.2) is valid for $N^*(h)$:

$$
\begin{aligned}
N^*(h) &= (nh^2)\iint \frac{\{\hat{p}^*(x,y) - \tilde{p}_{\tilde\theta^*}^*(x,y)\}^2}{R^2(K)p_{\tilde\theta}(x,y)}\omega(x,y)\,dx\,dy + \tilde{o}_p\{h\} \\
(\text{A.25}) \qquad &= N_1^*(h) + N_{2\tilde\theta}^*(h) + N_{3\tilde\theta}^*(h) + \tilde{o}_p(h),
\end{aligned}
$$

where $N_j^*(h)$ for $j = 1, 2$ and $3$ are the bootstrap versions of $N_j(h)$, respectively. As the bootstrap resample is generated according to $p_{\tilde\theta}$, the same arguments which led to Lemma 5 mean that $N_2^*(h) = N_3^*(h) = 0$. Thus, $N^*(h) = N_1^*(h) + \tilde{o}_p(h)$, where

$$
N_1^*(h) = (nh^2)\iint \frac{\{\hat{p}^*(x,y) - \tilde{p}^*(x,y)\}^2}{R^2(K)p_{\tilde\theta}(x,y)}\omega(x,y)\,dx\,dy.
$$

And similar lemmas to Lemmas A.6 and A.7 can be established to study $N_{1j}^*(h)$ which are the bootstrap versions of $N_{1j}^*(h)$, respectively.

PROOF OF THEOREM 2. It can be shown by taking the same route that establishes Lemma A.8 that, as $n \to \infty$ and conditioning on $\{X_t\}_{t=1}^{n+1}$, the distribution of $(L^*(h_1), \ldots, L^*(h_J))$ converges to $N_J(\beta\mathbf{1}_J, \Sigma_J)$ in probability, which readily imply the conclusion of the theorem. $\square$

Let $l_\alpha^*$ be the upper-$\alpha$ conditional quantile of $L_n^* = \max_{h\in\mathcal{H}} L^*(h)$ given $\{X_t\}_{t=1}^{n+1}$.

PROOF OF THEOREM 3. Let $l_{0\alpha}$ be the upper-$\alpha$ quantile of $\max_{1\le i\le J} Z_i$. From Theorem 2, and due to the use of the parametric bootstrap,

$$
(\text{A.26}) \qquad l_\alpha^* = l_{0\alpha} + o_p(1)
$$

under both $H_0$ and $H_1$. As $L_n \xrightarrow{d} \max_{1\le i\le J} Z_i$, by the Slutsky theorem,

$$
P(L_n \ge l_\alpha^*) = P(L_n + o_p(1) \ge l_{0\alpha}) \to P\left(\max_{1\le i\le J} Z_i \ge l_{0\alpha}\right) = \alpha,
$$

which completes the first part of Theorem 3. The second part of Theorem 3 is a direct consequence of Lemma A.10 and (A.26). $\square$

**Acknowledgments.** We thank Editor Jianqing Fan, an Associate Editor and two referees for their valuable comments and suggestions which have improved the presentation of the paper.

S. X. CHEN
C. Y. TANG
DEPARTMENT OF STATISTICS
IOWA STATE UNIVERSITY
AMES, IOWA 50011-1210
USA
E-MAIL: songchen@iastate.edu
yongtang@iastate.edu

J. GAO
SCHOOL OF MATHEMATICS AND STATISTICS
UNIVERSITY OF WESTERN AUSTRALIA
CRAWLEY, WESTERN AUSTRALIA 6009
AUSTRALIA
E-MAIL: jiti.gao@uwa.edu.au